\documentclass[preprint]{imsart}

\usepackage[T1]{fontenc}
\usepackage[babel]{csquotes}
\usepackage[utf8]{inputenc}
\usepackage[english]{babel}
\usepackage{amsfonts}
\usepackage{color}
\usepackage{amscd,amsthm}
\usepackage{mathrsfs,ae,dsfont,algpseudocode,algorithm}
\usepackage{graphicx}
\usepackage[normalem]{ulem}
\usepackage{amssymb,enumitem,fancyhdr}
\usepackage{nicefrac,ucs,xcolor,cmll}
\usepackage{mathrsfs,ae,dsfont,algpseudocode,algorithm}
  \expandafter\let\csname equation*\endcsname\relax
  \expandafter\let\csname endequation*\endcsname\relax
\usepackage{amsmath}
\usepackage{verbatim}
\usepackage{fullpage}
\definecolor{Red}{rgb}{1,0,0}
\definecolor{Blue}{rgb}{0,0,1}

\begin{document}
\begin{frontmatter}
\title{Bayesian smoothing of dipoles in Magneto--/Electro--encephalography}

\begin{aug}
\author{Valentina Vivaldi and Alberto Sorrentino}
\affiliation{Dipartimento di Matematica -- Universit\`a di Genova\\
Via Dodecaneso, $35$ -- $16146$ Genova (\@Italia\@) \\  \vspace{.2cm} and CNR--SPIN, Genova, Italy}

\ead{vivaldi@dima.unige.it, sorrentino@dima.unige.it}
\end{aug}
\begin{abstract}

We describe a novel method for dynamic estimation of multi--dipole states from Magneto/Electro--encephalography (M/EEG) time series.
The new approach builds on the recent development of particle filters for M/EEG; these algorithms approximate, with samples and weights, the posterior distribution of the neural sources at time $t$ given the data up to time $t$. However, for off--line inference purposes it is preferable 
to work with the smoothing distribution, i.e. the distribution for the neural sources at time $t$ conditioned on the whole time series.
In this study, we use a Monte Carlo algorithm to approximate the smoothing distribution for a time--varying set of current dipoles.
We show, using numerical simulations, that the estimates provided by the smoothing distribution are more accurate than those provided by the filtering distribution, particularly at the appearance of the source. We validate the proposed algorithm using an experimental dataset recorded from an epileptic patient. Improved localization of the source onset can be particularly relevant in source modeling of epileptic patients, where the source onset brings information on the epileptogenic zone.


\end{abstract}
\end{frontmatter}

\begin{comment}
If some form of the paper has appeared previously, even in a conference proceedings, this fact must be indicated clearly in both a cover letter and a footnote on the first page of the paper.\\

http://www.siam.org/journals/siap/authors.php
\end{comment}

\section{Introduction}
% versione alternativa:
%Non--invasive estimation of neural activity from external recordings, obtained by Magneto/Electro-Encephalography, are increasingly used in both neuroscience, for unraveling the mechanisms underlying cognitive processes, and in clinics, e.g., for localizing the epileptogenic areas. 

Magnetoencephalography and Electroencephalography (M/EEG) record non-invasively the magnetic field and the electric potential, respectively, produced by neural currents. M/EEG devices typically contain between few tens and few hundred of sensors, placed around the subject's head, that record the most direct consequence of the electrical brain activity at 1,000 Hertz. 
In the M/EEG inverse problem, one wants to find an estimate of the electrical current distribution inside the head. Thanks to their high temporal resolution, these estimates allow the investigation of the neural dynamics on a millisecond--per--millisecond scale, thus helping neuroscientists to discover how the brain performs higher cognitive functions, or helping clinicians to localize epileptogenic areas in patients. \\
In the \textit{dipolar} model, the neural generators of M/EEG data are described in terms of a set of point sources (named \textit{current dipoles} \cite{haetal93}): each dipole represents the activity of a small cortical area. Estimating the neural activity within the dipolar model requires to determine the number of dipoles, the dipole locations, orientations and strengths. As these parameters change relatively slowly in time (each dipole can remain active from few milliseconds up to several hundreds), there is benefit in using the a priori information that the neural currents change smoothly in time. Indeed, in the last decade there has been growing interest towards \textit{Bayesian filtering} \cite{sovoka03,mowisa08,soetal09,mietal13,soetal13,chsaso15}. Here the posterior distribution of the neural sources at time $t$ (the \textit{filtering} ditribution) is obtained from the posterior distribution of the neural sources at time $t-1$ with a two--step process. First, the prior distribution at time $t$ is obtained as the convolution of the posterior at time $t-1$ with a transition kernel that embodies a probabilistic evolution model for the neural sources, thus incorporating the smoothness prior; then, the posterior distribution at time $t$ is obtained via Bayes theorem. Due to the non--linearity of the forward model, these distributions cannot be calculated analytically; instead, they are sampled sequentially (in time) with a Monte Carlo importance sampling strategy, consisting in drawing a set of samples from the prior distribution and weighting them based on their likelihood. These algorithms, implementing Bayesian filtering with a Monte Carlo sampling technique \cite{roca04}, are usually referred to as \textit{particle filters}.\\
While the filtering distribution is particularly convenient for it is easy to approximate with Monte Carlo sampling, it is not the optimal choice for inference when on--line inference is not required. Indeed, the filtering distribution only embodies the information from the past time points, and completely ignores the information contained in the future time points. When possible, it is preferable to work with the \textit{smoothing} distribution, i.e. the distribution at time $t$ conditioned on the whole time series. 
It is foreseeable that the smoothing distribution will provide better estimates of brain activity particularly at the beginning of the activation, i.e. when the past contains little information on the source, while the future contains more. In this study, we build on previous work on particle filtering
and construct a Monte Carlo algorithm that approximates the smoothing distribution for a time--varying set of current dipoles.
\\
There are two well--known approaches to particle smoothing \cite{brdoma10}: forward filtering -- backward smoothing and the two--filter smoothing. The forward filtering -- backward smoothing consists in running a particle filter that approximates the filtering distribution, and then re-weight the samples going backward in time; in our context, the main limitation of this approach is that the approximation of the smoothing distribution makes use of the very same samples as that of the filtering distribution; therefore, if the filtering distribution has not found the high--probability region, the smoothing distribution will also be poorly approximated.
The two--filter smoothing runs two separate filters, that approximate the filtering distribution and the backward information filter, and then re-weights the samples from the backward information filter to approximate the smoothing distribution. Because it uses only the samples of the backward filtering, it has the usual disadvantages: if the backward filtering has only found low-probability regions the smoothing distribution will be poorly approximated.
In addition, both approaches suffer from a computational cost scaling quadratically with the number of particles.

In this paper we propose a slightly modified two--filter smoother, that tries to overcome the limitations previously described. The idea is that one can modify the recursion of the two--filter smoothing, and obtain a second approximation of the smoothing distribution which makes use of the samples of the forward filter.
This way, at each time we obtain two (possibly different) approximations of the smoothing distribution; we then select the best approximation based on the marginal likelihood of the underlying filtering density. Furthermore, in order to reduce the overall computational cost, we run the forward and backward filters with a large number of particles, but then we sub--sample these distributions retaining only a small fraction of the sample set.
As a result, the algorithm keeps a reasonable computational cost, seems to be more effective in detecting the dipole sources at their onset, and can find useful applications with epileptic patients data in which estimating the onset is crucial for identifying the epileptogenic zone.\\

The paper is organised as follows: in Section \ref{Sec:PFS} we review the Bayesian dynamic dipole models for Magnetoencephalography; in Section \ref{Sec:DSA}, we introduce the smoothing algorithms and we describe our double two--filter smoothing; in Section \ref{Sec:DSA} we describe how to apply the double two--filter smoothing algorithm to the M/EEG problem; Section \ref{Sec:SIM} and \ref{Sec:REA} provide a validation of our algorithm via a simulation study and an illustration of performance on real data. Our conclusions are offered in Section \ref{Sec:DIS}.

\section{Bayesian dynamic dipole models for Magnetoencephalography}
\label{Sec:PFS}

\subsection{Magneto--/Electro--encephalography}

The definition of the forward model for the M/EEG signals is based on the quasi--static approximation \cite{sa87} of Maxwell equations, in which the electrical currents produced by the neuronal discharges play the role of the source term. A detailed treatment of the physics can be found e.g. in \cite{haetal93} and \cite{sa87}. 
Here we just describe the mathematical model we adopt for the source term and for the forward problem.

We model the neural sources as the superposition of an unknown number $N$ of current dipoles \cite{soluar14}. Each current dipole is parameterized by a location $r$ in the brain and a dipole moment $q$, representing orientation and strength of the current at location $r$. 
For practical and computational reasons, source locations belong to a pre--defined grid $R_{grid}$; dipole moments are three--component vectors. The state--space for a single dipole is therefore $\mathcal{D} = R_{grid} \times \mathbb{R}^3$. The state--space $\mathcal{D}_N$ for a fixed number $N$ of dipoles can be obtained as the Cartesian product of $N$ single--dipole state spaces, $\mathcal{D}_N = \mathcal{D} \times \dots \times \mathcal{D}$. Since the number of dipoles is unknown, the state space $\mathcal{J}$ of the neural current is constructed by union of spaces with a fixed number of dipoles, i.e.

\begin{equation}
\mathcal{J}:= \bigcup_{N=0}^{N_{max}}\{N\} \times \mathcal{D}_N \backslash \sim
\end{equation}
where $\mathcal{D}_0 = \emptyset$, and $\sim$ is an equivalence relation that accounts for the fact that two states that only differ by a permutation of the ordering of the dipoles are physically equivalent. 
In a M/EEG time series, one can either model the dipoles as stationary (dipole locations don't change in time, while dipole moments do) or moving (dipole locations change in time, as well as dipole moments). The case of stationary dipoles has been recently addressed in \cite{soetal13,soso14}; here we consider the case of moving dipoles, i.e. the neural current at time $t$ is:
\begin{equation}
j_t = ( N_t, \left\{ r_t^i, q_t^i \right\}_{i=1,...,N_t}).
\label{eq:dipoles}
\end{equation}

\noindent
The data $d_t$ produced by a set of dipoles such as (\ref{eq:dipoles}) is the superposition of the data produced by individual dipoles, i.e.: 
\begin{equation}\label{eq:forward_model}
d_t = \sum_{i=1}^{N_t}G(r_t^i) \cdot q_t^i + \epsilon_t
\end{equation}
where $\epsilon_t \sim \mathcal{N}(0,\sigma_{noise})$ is the noise component that is assumed to be additive, zero-mean and Gaussian; $G(r)$ is the \textit{leadfield} associated to location $r$, and can be thought of as a 3--column matrix containing the data produced by a unit dipole located in $r$ and oriented along the three orthogonal directions. It is common practice to pre--compute $G(r)$ for a large set of points (typically around 10,000 points) distributed inside the brain volume and store the result in a large matrix; therefore, source locations are practically constrained on a grid.

\subsection{Bayesian particle filtering}

The inverse source problem of M/EEG consists in the estimation of the neural current, knowing the measured data. The literature concerning this problem is extensive \cite{uuhasa98,soka04,scgewo99}. We have considered a Bayesian approach in which all variables are modeled as random variables.\\
In a Bayesian setting, the inverse M/EEG problem can be casted as an ill--posed dynamic inverse problem, through a Hidden Markov Model \cite{camory05}.
Indeed, one can model the neural current  $j_t$ and the data $d_t$ as two Markov processes ${\{J_t\}}_{t=1}^T$ and $\{D_t\}_{t=1}^T$ satisfying

\begin{eqnarray}
p(j_{t+1}|j_t, j_{t-1},...,j_1, ) = p(j_{t+1}|j_t) \label{transition}\\
p(j_{t+1}|j_t, d_t,...,d_1) = p(j_{t+1}|j_t)\\
p(d_t|j_t, j_{t-1},...,j_1, d_{t-1},...,d_1) = p(d_t|j_t) \label{likelihood}
\end{eqnarray}
i.e., the neural current are a first--order Markov process, and the measured data are a first--order Markov process with respect to the history of $J$; we also assume that the two processes are homogeneous, i.e. the distributions (\ref{transition}) and (\ref{likelihood}) do not change in time. Then, the problem of \textit{filtering} is the one to infer information about the state of the current at time $t$, given the data up to time $t$,
i.e. to obtain the filtering distribution $p(j_t|d_{1:t})$. Given a prior distribution at the first time sample $p(j_1)$, the transition kernel $p(j_{t+1}|j_t)$ and the likelihood $p(d_t|j_t)$, the problem
can be solved by sequential application of a two--step algorithm, known as \emph{Bayesian filtering}; the posterior distribution at time $t$
is computed by Bayes theorem
\begin{equation}\label{eq:bayes_filtering}
p(j_t|d_{1:t}) = \frac{p(d_t|j_t) p(j_t|d_{1:t-1})}{ \int p(d_t|j_t)p(j_t|d_{1:t-1})\, dj_t}
\end{equation}
where $d_{1:t}:=(d_1, \dots d_t)$, while $p(j_t|d_{1:t-1})$ plays the role of prior at time $t$ and $p(j_1|d_{1:0}) := p(j_1)$. The next prior
is computed by means of the Chapman--Kolmogorov equation

\begin{equation}\label{eq:bayes_filtering2}
p(j_{t+1}|d_{1:t}) = \int p(j_t|d_{1:t}) p(j_{t+1}|j_t) dj_t~~~.
\end{equation}
For linear--Gaussian models, these two formulas lead to the well--known Kalman filter \cite{ka60}. For non--linear/non--Gaussian models, numerical approximations, like particle filters, are needed.\\

\noindent
Particle filters are receiving growing attention in the last years \cite{sovoka03, dojo11,paetal07,sopapi07}. They are a Monte Carlo technique that rely on a sequential application of an \textit{importance sampling} scheme: the filtering distribution is approximated with a weighted set of samples
\begin{equation}
p(j_t|d_{1:t}) \simeq \sum_{l=1}^\alpha w_t^l \delta(j_t,j_t^l)
\end{equation}
where $j_t^l$ are the sample points (particles), $w_t^l$ are the weights and $\delta(\cdot, \cdot)$ is the Kronecker delta. 
Samples $j_{t}^i$ are drawn from an importance distribution $\eta(j_t|d_{1:t-1})$; then the weights for approximating the posterior are given by the ratio $w_t^l = \frac{p(j_t^l|d_{1:t})}{\eta(j_t^l|d_{1:t-1})}$; when these weights turn out to be too diverse, so that most of them are negligible, one can resample this weighted set by taking multiple times particle with high weights, and discarding particles with low weigths. One of the simplest implementations consists in choosing the importance distribution equal to the marginal prior distribution $\eta(j_t|d_{1:t-1}) = p(j_t|d_{1:t-1})$; this case is often referred to as SIR (Sampling Importance Resampling) particle filter; in the simulations below we will be using a slightly modified version of SIR \cite{dogoan00}.

\section{Particle Smoothing of Dipoles}

\label{Sec:DSA}

The smoothing problem is the one to make inference about the state at time $t$, given the whole sequence of measurements, up to the final time point $T$; i.e., one is interested in the distributions $p(j_{t}|d_{1:T})$. 
A common approach, usually referred to as two--filter smoothing, consists in exploiting the following identity
\begin{equation}\label{eq:oursmoothing}
p(j_t|d_{1:T}) = {p(j_t|d_{1:t-1}) p(d_{t:T}|j_t) \over p(d_{t:T}|d_{1:t-1})}~~~.
\end{equation}
However, Monte Carlo approximation of these densities is not straightforward. Indeed, the first term at the numerator of the right hand side is routinely approximated by a particle filter; but the second term, often referred to as \emph{backward information filter}, is not a probability density with respect to $j_t$. In the following subsection we summarize, for the sake of clarity, the approach proposed in \cite{brdoma10}.

\subsection{Two--filter smoothing}

In order to approximate the backward filter, in \cite{brdoma10} the authors introduce a set of auxiliary densities $\gamma_t(j_t)$ and probability distributions $\tilde{p}(j_t|d_{t:T})$ such that 

\begin{equation}\label{eq:smoothing_rel_prop}
\tilde{p}(j_t|d_{t:T}) \propto p(d_{t:T}|j_t) \gamma_t(j_t)~~~;
\end{equation}

so that we can re--write
\begin{equation}\label{eq:smoothing_a}
p(j_t|d_{1:T}) \propto {p(j_t|d_{1:t-1}) \tilde{p}(j_t|d_{t:T}) \over \gamma_t(j_t)}~~~.
\end{equation}

From equation (\ref{eq:smoothing_rel_prop}), one can derive a recursion for $\tilde{p}(j_t|d_{t:T})$ that can be conveniently written as a two--step algorithm that closely reminds the filtering equations (\ref{eq:bayes_filtering})-(\ref{eq:bayes_filtering2}):

\begin{equation}\label{eq:s1}
\tilde{p}(j_t | d_{t+1:T}) = \int \tilde{p}(j_{t+1}|d_{t+1:T}) \frac{p(j_{t+1}|j_t) \gamma_t(j_t)}{\gamma_{t+1}(j_{t+1})} dj_{t+1}
\end{equation}

\begin{equation}\label{eq:s2}
\tilde{p}(j_t | d_{t:T}) = \frac{p(d_t|j_t) \tilde{p}(j_t | d_{t+1:T})}{\int p(d_t|j_t) \tilde{p}(j_t | d_{t+1:T}) dj_t} ~~~.
\end{equation}

These equations allow to use a particle filter, going backward in time, to approximate the backward information filter with a set of samples $\tilde{j}_t^i$ and weights $\tilde{w}_t^i$.

After having obtained an approximation for the filtering density and for the backward information, one can plug both in (\ref{eq:smoothing_a}) and use (\ref{eq:bayes_filtering2}) to obtain:

\begin{eqnarray}
	p(j_t|d_{1:T}) &\propto& {p(j_t|d_{1:t-1}) \tilde{p}(j_t|d_{t:T}) \over \gamma_t(j_t)} = \nonumber \\
	&=& {\int p(j_t|j_{t-1}) p(j_{t-1}|d_{1:t-1}) \tilde{p}(j_t|d_{t:T}) dj_{t-1} \over \gamma_t(j_t)} ~~~. 
\end{eqnarray}
Replacing by their weigthed samples $p(j_{t-1}|d_{1:t-1}) \simeq \sum_{k=1}^\alpha w_{t-1}^k \delta(j_{t-1}, j_{t-1}^k)$ and $\tilde {p}(j_t|d_{t:T}) \simeq \sum_{l=1}^\alpha \tilde{w}_t^l \delta(j_t, \tilde{j}_t^l) $ , one obtains the approximation of the smoothing distribution provided by the two--filter smoothing:
\begin{equation}\label{eq:sm}
p_{1}(j_t|d_{1:T}) = \sum_{l=1}^\alpha w_{1,t|T}^l \delta(j_t, \tilde{j_t^l})
\end{equation}
where
\begin{equation}
w_{1,t|T}^l \propto \tilde{w}_t^l \left( \sum_{k=1}^\alpha w_{t-1}^k {p(\tilde{j}_t^l|j_{t-1}^k) \over \gamma_t(\tilde{j}_t^l)} \right)
\label{eq:s_weight}
\end{equation}
Namely, the two--filter smoothing uses the same particles used to approximate the backward filter; importantly, the computation of the new weights (\ref{eq:s_weight}) requires $\alpha^2$ operations.

\subsection{A double two--filter smoothing}

The algorithm we propose is a modified version of the two--filter smoothing; first, we observe that it is possible to re--write (\ref{eq:oursmoothing}) by moving $d_t$ from the second to the first term at the right hand side:
\begin{equation}\label{eq:oursmoothing_bis}
p(j_t|d_{1:T}) = {p(j_t|d_{1:t}) p(d_{t+1:T}|j_t) \over p(d_{t+1:T}|d_{1:t})} ~~~;
\end{equation}
then,  by using (\ref{eq:smoothing_rel_prop}) and (\ref{eq:s2}) we obtain the following relation 
\begin{equation}\label{eq:smoothing_a2}
p(j_t|d_{1:T}) \propto { p(j_t|d_{1:t}) \tilde{p}(j_t|d_{t+1:T}) \over \tilde{p}_0(j_t) }
\end{equation}
Like in the two--filter smoothing, $p(j_t|d_{1:t})$ and $\tilde{p}(j_t|d_{t+1:T})$ can be approximated by two particle filters, one going forward and one backward in time; then we can use (\ref{eq:s1}) to obtain the following approximation of the smoothing distribution:

\begin{equation}\label{eq:sm_2}
p_{2}(j_t|d_{1:T}) = \sum_{l=1}^\alpha w_{2, t|T}^l \delta(j_t, j_t^l)
\end{equation}
where
\begin{equation}
w_{2,t|T}^l \propto w_t^l \left( \sum_{k=1}^\alpha \tilde{w}_{t+1}^k  {p(\tilde{j}_{t+1}^k|j_t^l) \over \gamma_{t+1}(\tilde{j}_{t+1}^k)} \right)
\label{eq:s_weight_2}
\end{equation}

namely, this time the smoothing distribution is approximated using the same particles used to approximate the forward filter; the structure of the weights in (\ref{eq:sm_2}) is analogous to that of (\ref{eq:sm}), i.e. the computational cost is again $\alpha^2$.\\

\noindent
From time to time, the particles coming either from the forward or the backward filter are not well suited for approximating the smoothing distribution. Here we partially overcome this problem by proposing an algorithm, we call it \emph{double two--filter smoothing}, that uses both approximations (\ref{eq:sm}) and (\ref{eq:sm_2}). Importantly, we do not combine the two samples, which is difficult because the normalizing constants are unknown, but we limit ourselves to selecting one of the two approximations; we reckon there may be several ways to select which approximation is to be preferred, therefore we describe the criterion we use in the following section. The resulting algorithm runs as follows:

\begin{enumerate}
	\item run the forward and the backward filters;
	\item approximate the smoothing distributions with (\ref{eq:sm}) and (\ref{eq:sm_2});
	\item for each time $t$ pick the approximation that best explains the data and set either $p_{dtfs}(j_t|d_{1:T}) = p_{1}(j_t|d_{1:T})$ or $p_{dtfs}(j_t|d_{1:T}) = p_{2}(j_t|d_{1:T})$.
\end{enumerate}

\section{Application to M/EEG}

In order to apply the double two--filter smoothing to the MEEG problem, we need to devise: (i) the statistical model, i.e. the prior distribution, the transition kernel and the likelihood function; (ii) the algorithm settings, i.e. the importance densities for the forward and backward filters, the auxiliary densities and a criterion to select the best approximation.

\subsection{Statistical Model} 

\textbf{Initial prior distribution}. We set the initial prior distribution based on neurophysiological considerations. In general, the number of active dipoles is expected to be small, i.e. between 1 and 5 ($=N_{max}$); therefore we use a Poisson prior for $N_1$ with rate parameter below 1. Conditional on the number of dipoles, the dipole parameters are independent. The prior for the dipole locations is uniform in the brain volume; we recall that for computational reasons dipole locations are constrained to a finite set of values. The prior for the dipole orientation is uniform in the sphere, and the prior for the dipole strength is log--uniform. As a result, the initial prior distribution can be written as
\begin{equation}\label{eq:prior}
p(j_1) = \sum_{n=0}^{N_{max}} \mathbb{P}(N_1 = n) \prod_{m=1}^n \mathcal{U}_{R_{grid}}(r_1^{m}) \mathcal{U}_S\left(\frac{q_1^m}{|q_1^m|}\right) \mathcal{LU}(|q_1^m|) .
\end{equation}
where $\mathcal{U}_{R_{grid}}$ is the uniform distribution over the grid $R_{grid}$, discretizing the brain volume; $\mathcal{U}_S$ is the uniform prior on the spherical surface; $\mathcal{LU}$ is the log--uniform prior for the dipole strength.\\

\noindent
\textbf{Transition kernel}. In our model, at each time point a new dipole can appear, and existing dipoles may disappear; in addition, dipole locations, orientations and strengths can change. Therefore the transition density accounts for the possibility of dipole birth, dipole death and evolution of the dipole parameters. 
To limit the complexity of the model, only one birth and one death can happen at any time point; due to the high temporal resolution of the recordings, this is not a real limitation, i.e., several dipoles can appear and disappear in few milliseconds. Our transition density can be written as follows:
\begin{multline}\label{eq:TK}
p(j_{t+1} | j_{t}) =    \\
P_{\rm birth} \times U_{R_{\rm grid}}(r_{t+1}^{N_{t+1}}) \mathcal{U}_S\left(q_{t+1}^{N_{t+1}} \over |q_{t+1}^{N_{t+1}}| \right) \mathcal{LU}(q_{t+1}^{N_{t+1}}) \times \prod_{i=1}^{N_{t}} M(r_{t+1}^{i},r_{t}^{i}) \mathcal{N} (q_{t+1}^{i}; q^{i}_{t},\Delta^{i}_{t} ) + \\
+ P_{\rm death} \times \frac{1}{N_{t}} \sum_{j=1}^{N_{t}} \prod_{i=1}^{N_{t}-1} M(r_{t+1}^{i},r_{t}^{a_{j,i}}) \mathcal{N} (q_{t+1}^{i}; q_{t}^{a_{j,i}}, \Delta^{a_{j,i}}_{t})  + \\		
+ (1-P_{\rm birth} - P_{\rm death}) \times \prod_{i=1}^{N_{t}} M(r_{t+1}^{i},r_{t}^{i}) \mathcal{N}(q_{t+1}^{i};q^{i}_{t},\Delta^{i}_{t} )\,.				
\end{multline}
where $M(r_1,r_2)$ represents the transition probability from location $r_1$ to location $r_2$; in the simulations below we use
$$M(r_1,r_2) \propto e^{- {|| r_1-r_2||^2 \over 2 \cdot \rho^2}}$$ 
i.e., the transition probability is proportional to a Gaussian centered at the current location; the value of the standard deviation $\rho$ is set to 5 mm, which worked fine with our $10^4$-points discretization of the brain volume.\\
The first term in (\ref{eq:TK}) takes into account the chance that a new dipole appears, with probability $P_{birth}$; the new dipole location is then uniform in the grid, the orientation is uniform in the sphere, and the strength is log--uniform. All other dipoles evolve independently: dipole locations change according to matrix $M$, while dipole moments perform a Gaussian random walk in which $\Delta^{i}_{t}$ has been set to ${ ||q_{t}^{i}||_{2}\over 5}$, where $|| \cdot ||_2$ is Euclidean norm; such variable standard deviation allows the dipole strength to change non--negligibly but not too much in between two time points. \\
\noindent
The second term accounts for the possibility that one of the existing dipoles disappears: all the dipoles have the same possibility to disappear and the disappearance of a dipole entails a re-arrangement of the dipole labels that is given by
\begin{equation}
a_{j,m}=
\begin{cases}
m &\text{se $m$ < j}\\
m+1 & \text{se $m$ $\geq$ j.}
\end{cases}
\end{equation}
Finally in the last term the number of dipoles in the set remains the same.
Birth and death probabilities were set to $P_{birth} = 1/100$ and $P_{death} = (1-{(1-1/30)}^{N_t})$ respectively, as the expected lifetime of a single dipole is about 30 time points since simultaneous death are neglected.\\

\noindent
\textbf{Likelihood function}. Noise is assumed to be zero--mean Gaussian and additive. Therefore, the likelihood function is 

\begin{equation}
p(d_t|j_t) = \mathcal{N}\left( d_t; \sum_{i=1}^{N_t} G(r_t^i) \cdot q_t^i, \Sigma\right)
\end{equation}
where $\Sigma$ is the noise covariance matrix.
\subsection{Algorithm settings}

\textbf{Importance density for the forward filter}.
The choice of the importance distribution is known to play an important role in making a particle filter efficient. The simplest particle filter consists of using the marginal prior (\ref{eq:bayes_filtering2}) as importance distribution, but this is known to be not the optimal choice. In order to improve against this simple choice, while maintaining the same computational cost, we choose $\eta(j_t|d_{1:t-1})$ to be a modified version of (\ref{eq:bayes_filtering2}), where the transition kernel $p(\cdot | \cdot)$ is replaced by a kernel $\eta(\cdot | \cdot)$ with the same analytical form, but having $Q_{birth} = 1/3$ and $Q_{death}=1/3$ instead of $P_{birth}$ and $P_{death}$:

\begin{equation}
	\eta(j_t|d_{1:t-1}) = \int \eta(j_t|j_{t-1}) p(j_{t-1}|d_{1:t-1})  dj_{t-1}
\end{equation}

 This choice allows a better trans--dimensional sampling, which is particularly useful to explore the state--space of new--born dipoles.\\

\noindent
\textbf{Auxiliary and importance densities for the backward filter}.
For the second Monte Carlo filter that goes backward in time, we have to select the auxiliary distributions $\gamma_t(j_t)$ and the importance distributions $\eta(j_t|d_{t+1:T})$. For simplicity, we have chosen the auxiliary distributions to be all equal to the initial prior distribution, i.e. $\gamma_t(j_t) = p(j_1)$ for $t=1,...,T$. 
The importance distribution, on the other hand, has been set in analogy with that of the forward filter
to be

\begin{equation}
	\eta(j_t|d_{t+1:T}) = \int \eta(j_t|j_{t+1}) \tilde{p}(j_{t+1}|d_{t+1:T})  dj_{t+1}
\end{equation}
where $\eta(\cdot|\cdot)$, this time, goes backward.\\

\noindent
\textbf{Picking the smoothing distribution}. Once the forward and backward filters have been approximated, we compute two separate approximations of the smoothing distributions using (\ref{eq:sm}) and (\ref{eq:sm_2}). The last step of the double two--filter smoothing consists of choosing either of the two according to some criterion. 
Here, we select the approximation based on the marginal likelihood of the underlying filtering algorithm: indeed, we compute $L^f_t = \sum_{l=1}^{\alpha} W_t^l$ for the forward filter and $L^b_t = \sum_{l=1}^{\alpha} \tilde{W}_t^l$ for the backward filter, where $W_t^l$ and $\tilde{W}_t^l$ are the un--normalised particle weigths; we then pick the smoothing distribution that is based on the particle set that has obtained higher marginal likelihood in the filtering.\\
\noindent
To reduce the high computational cost associated with the calculation of the weights in eq. (\ref{eq:s_weight}) - (\ref{eq:s_weight_2}), in which the number of operations is proportional to $\alpha^2$, we choose a subset of only 100 samples from the whole particle set. Particles of this subset are sampled from the multinomial distribution defined by the forward and backward filtering weights; we note that, as a consequence of this choice, the weights $w_t$ and $\tilde{w}_t$ in the general formulas (\ref{eq:s_weight}) and (\ref{eq:s_weight_2}) are uniform.  We have observed that this subsampling does not affect the effectiveness of the algorithm; this is likely due to the relatively low number of particles with non--negligible weight in the (full--sample) forward and backward filtering. Furthermore, our approach might be seen as a computationally inexpensive and easy--to--implement approximation to the approaches described in \cite{kletal06} for reducing the computational cost associated with the $\alpha^2$ operations.

\section{Simulations}\label{Sec:SIM}
In this section MEG simulated data are used to validate the performance of the double two--filter smoothing approach. Here we describe the generation of synthetic data and we show our results.\\

\subsection{Generation of synthetic data}

Generation of synthetic time series was performed according to the following general scheme:
\begin{enumerate}
	\item a head model is defined, accounting for the geometrical and physical properties of the subject's head; we used a healthy subject's MRI to define a realistic head model; Freesurfer and MNE (http://www.martinos.org/mne/) were used to obtain a tessellation of the cortical surface and to compute the Boundary Element Method needed to obtain the leadfield matrix $G(r)$ (see eq. (\ref{eq:forward_model})) for $r \in R_{grid}$; the brain volume was discretized with 12,324 points;
	\item the spatio--temporal evolution of the neural sources is simulated, i.e. the variables $N_t$ and $\{r_t^i, q_t^i\}_{i=1,...,N_t}$ are assigned values; the generation of these values followed different criterions in Simulations (1)-(8), these criterions are described below;
	\item synthetic MEG recordings are generated according to equation (\ref{eq:forward_model}): first the noise--free data are generated through the leadfield matrix, then white Gaussian noise of fixed standard deviation is added to the noise--free data; while the noise standard deviation was fixed for all the simulations, different synthetic data sets can have remarkably different signal--to--noise ratios, because the strength of the noise--free signal can be very different.
\end{enumerate}

To validate our algorithm using a variety of synthetic experimental conditions, we devised eight different groups of Simulations, differing from each other for the number of sources (one or two), the dynamics of the source location (fixed or moving) and the dynamics of the dipole moment (fixed or moving). Each group of Simulations contains twenty simulations; in each simulation, the source locations and orientations are randomly drawn in the brain volume at $t=1$, and then evolve differently depending on the Simulation group. All simulations are 30 time points long. 

\begin{itemize}
   \item Simulation 1: one source with fixed location and fixed dipole moment;
   \item Simulation 2: one source with moving location and fixed dipole moment;
   \item Simulation 3: two sources with fixed location and fixed dipole moment;
   \item Simulation 4: two sources with moving location and fixed dipole moment;
   \item Simulation 5: one source with fixed location and bell--shaped dipole moment;
   \item Simulation 6: one source with moving location and bell--shaped dipole moment;
   \item Simulation 7: two sources with fixed location and bell--shaped dipole moment;
   \item Simulation 8: two sources with moving location and bell--shaped dipole moment;
\end{itemize}

\begin{figure}[ht!]
\centering
\includegraphics[width=17cm]{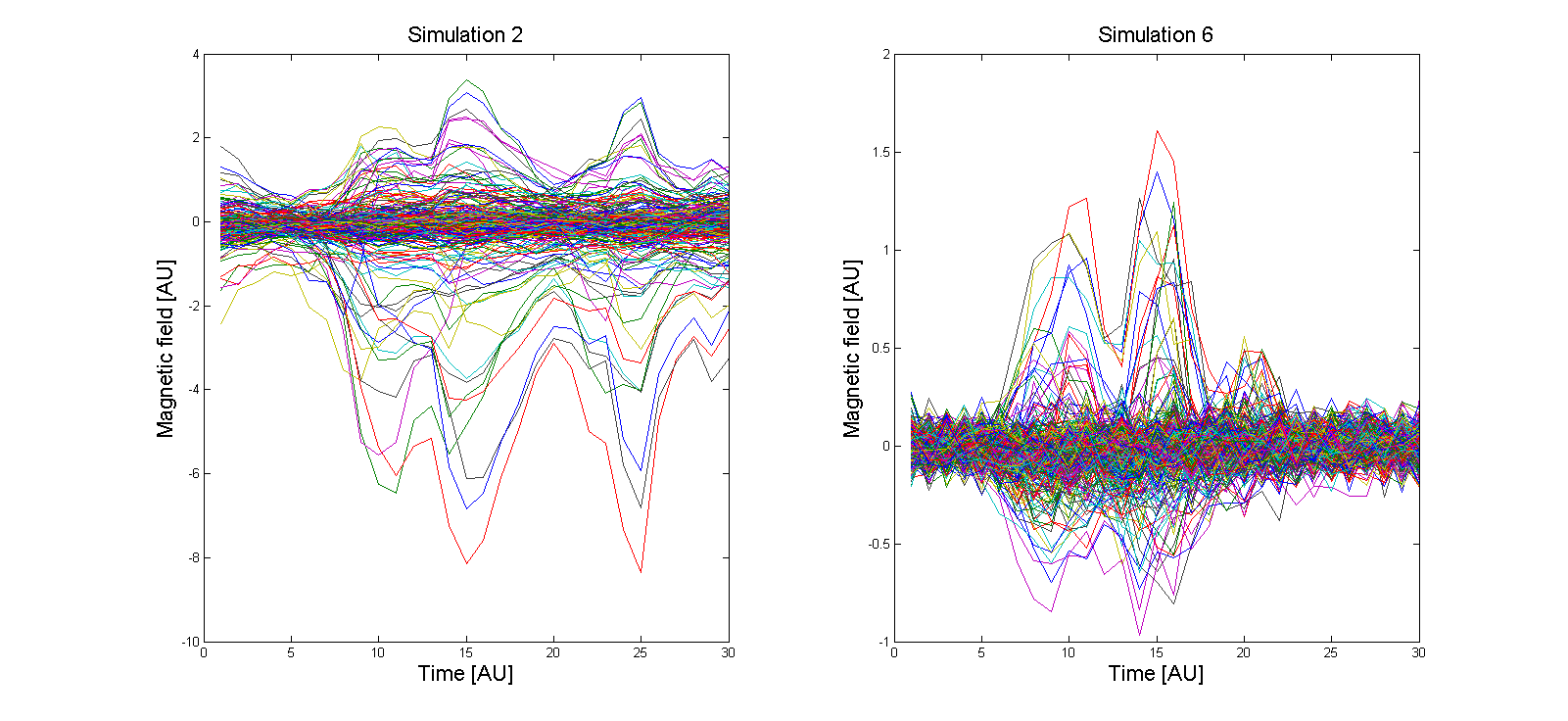}
\caption{Two synthetic time series: on the left panel, a data set generated by a single dipole with moving location and fixed dipole moment (Simulation 2); on the right panel, a data set produced by one dipole with moving location and bell--shaped dipole moment (Simulation 6).}
\label{fig_mag_12}
\end{figure}

Moving dipoles are generated by a random walk in the brain volume: the dipole location at time $t+1$ is constrained to the grid points that belong to the ball centered in the current location and of small radius (1 cm); in addition, the dipole location $r_{t+1}$ must satisfy $|r_{t+1}-r_{t-1}| > |r_t - r_{t-1}|$, so that dipoles cannot oscillate around the same position.
Bell--shaped dipole moments are generated using a Gaussian function, i.e. we set $|q_t| \propto \exp\{-(t-t_0)^2/2a^2\}$ with $t_0=15$ and $a=4$.
In Figure \ref{fig_mag_12} we show two examples of synthetic data used in our simulations.

\subsection{Point Estimates}

In order to evaluate the performance of the proposed method, point estimates are computed from the approximations to the filtering and smoothing distributions as follows.
Let $\{j_t^i, w_t^i\}$ be a sample approximating the smoothing distribution; following recent literature \cite{soetal13}, point estimates of the neural current parameters are computed as follows:
\begin{itemize}
	\item an estimate of the number of sources is obtained as the mode of the distribution 
	\begin{equation}
		\mathbb{P}(N_t=k|d_{1:T}) = \sum_{l=1}^\alpha w_t^l \delta(k,N_t^l)
	\end{equation}
	\item estimates of the source locations are obtained as the peaks (modes) of the intensity measure, defined as: 
	\begin{equation}
		p_t(r|d_{1:T})= \sum_{l=1}^\alpha w_t^{l} \sum_{k=1}^{N_t^l} \delta(r, r_t^{l, (k)})
	\end{equation}

	\item estimates of the dipole moments are obtained as the mean values of the conditional distributions
	\begin{equation}
		\mathbb{E}[q_t|r] = \sum_{l=1}^\alpha w_t^l \sum_{k=1}^{N_t^l} q_t^{l,(k)} \delta(r, r_t^{l,(k)})
	\end{equation}
\end{itemize}

\subsection{Results}

\begin{center}
\begin{figure}[H]
\includegraphics[width=8cm]{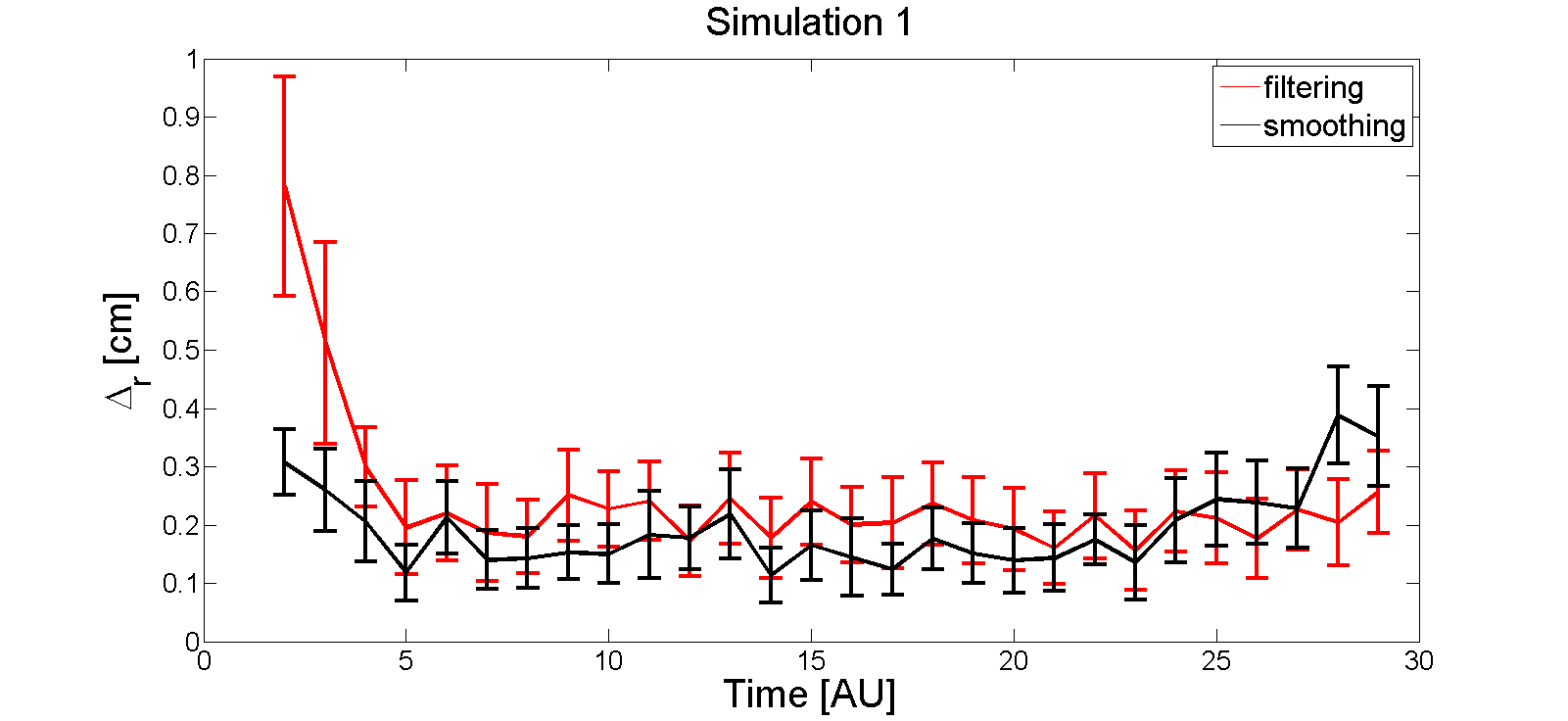}
\includegraphics[width=8cm]{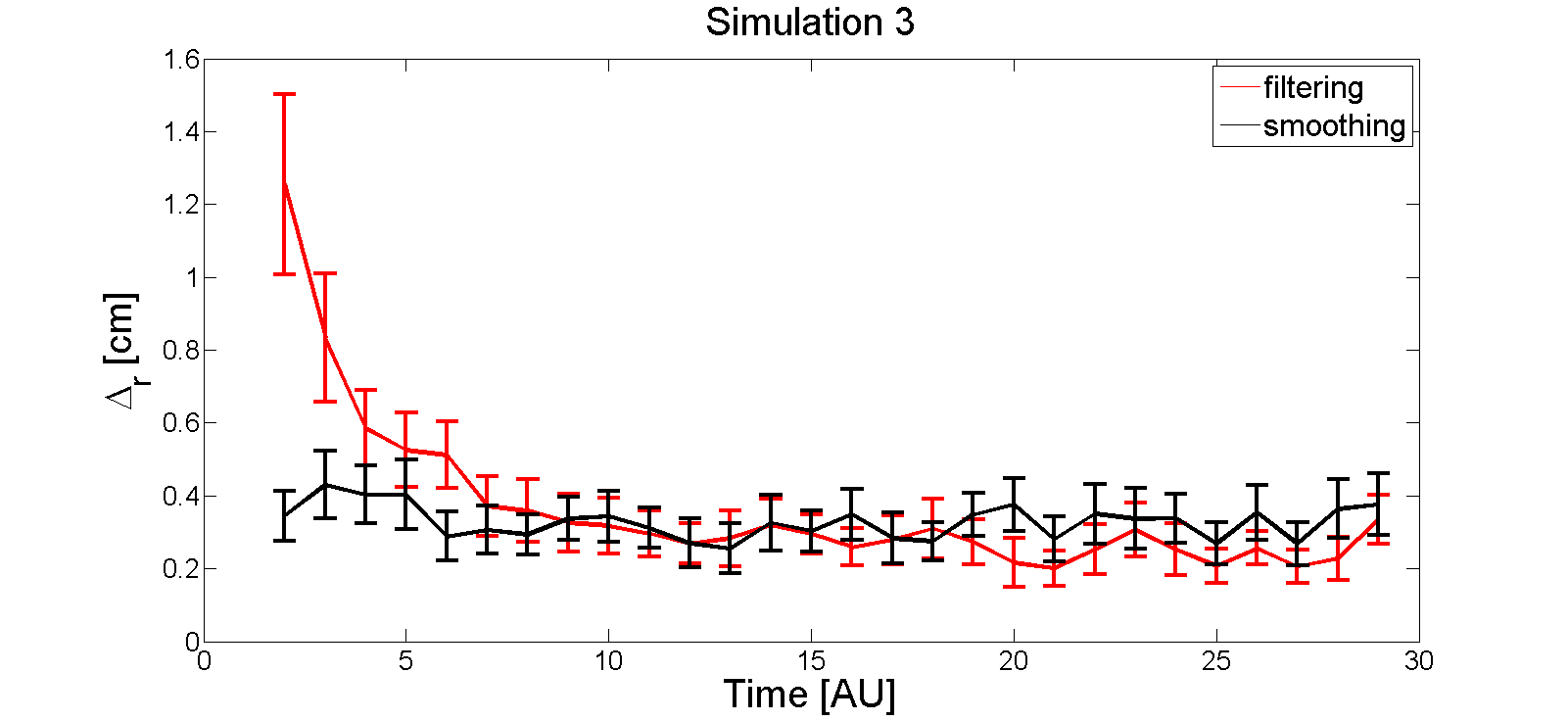}\\ 
\includegraphics[width=8cm]{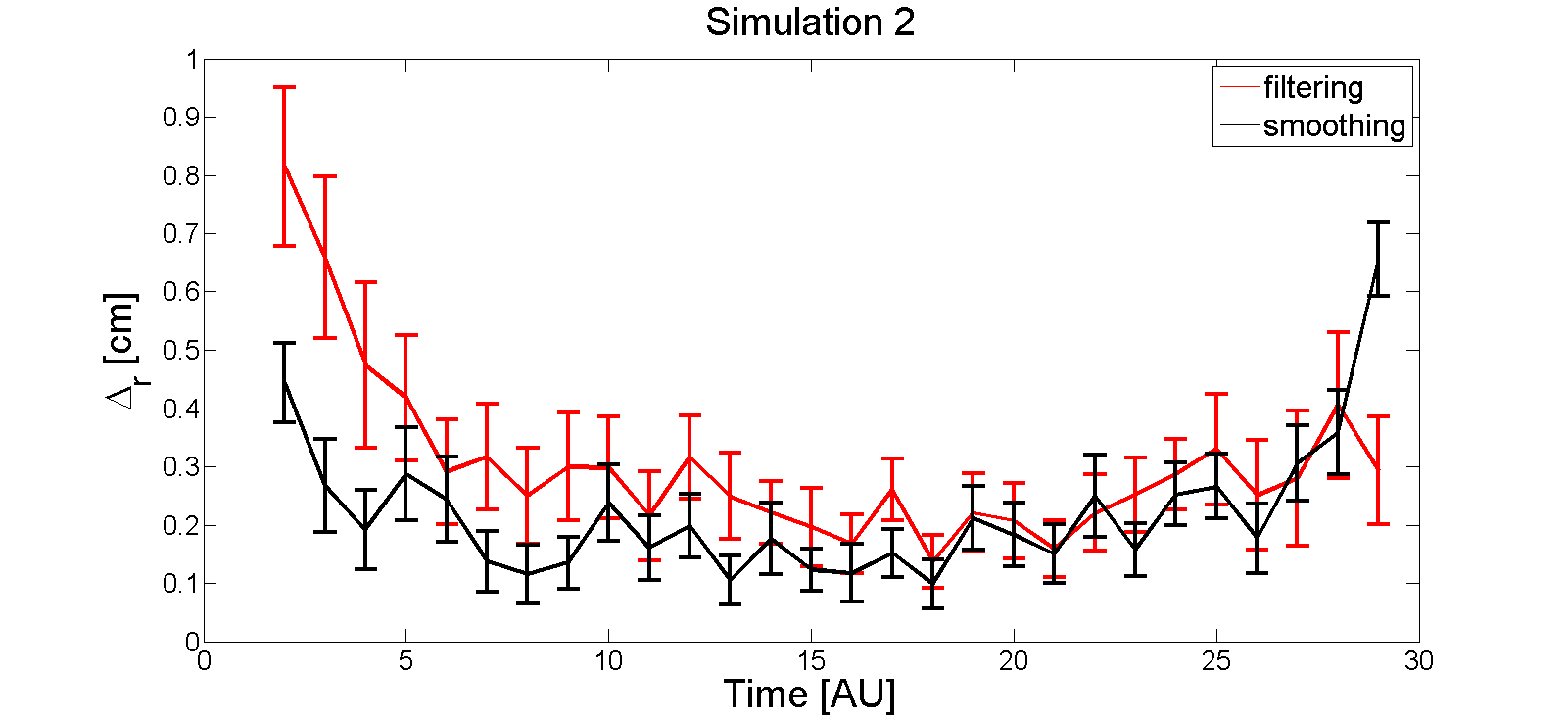}
\includegraphics[width=8cm]{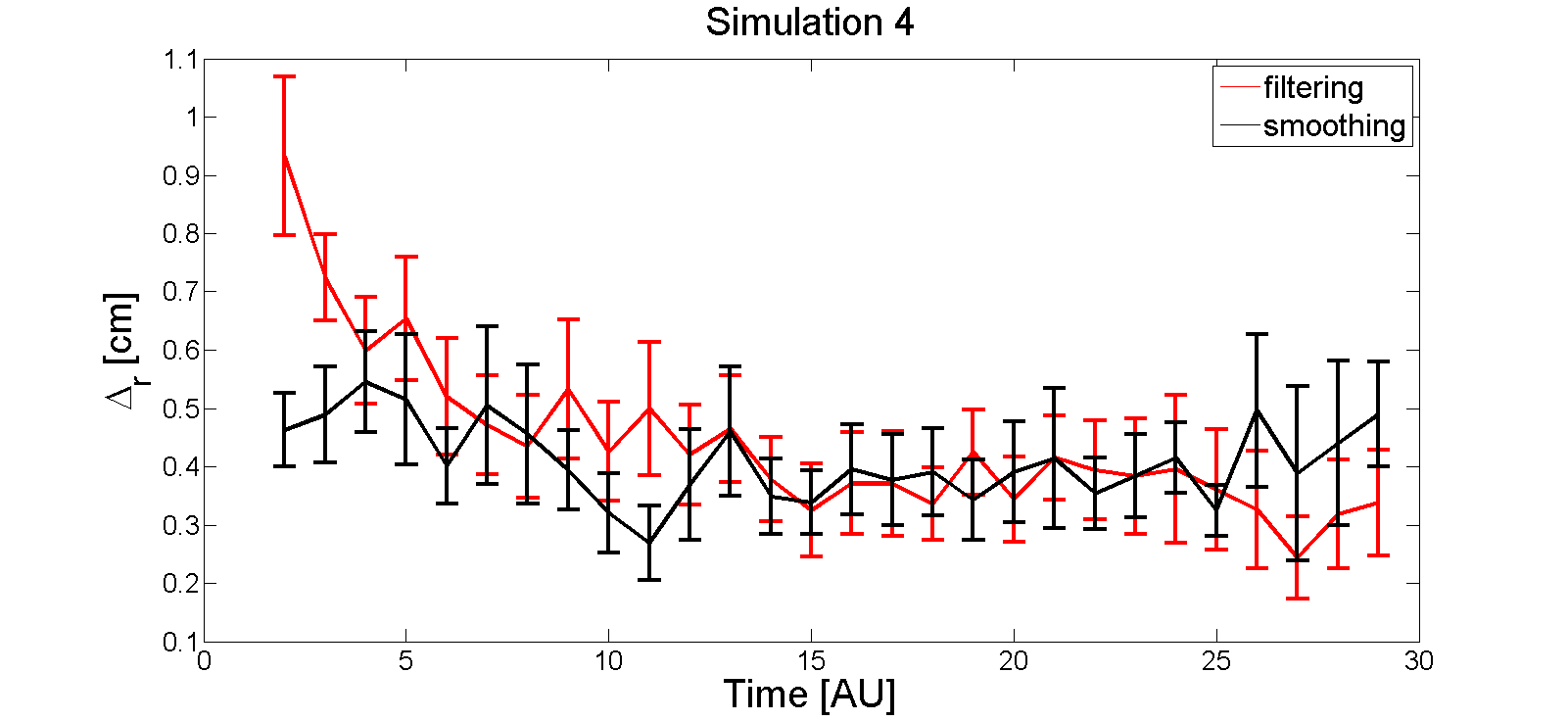}
\caption{Simulation 1--4: Mean localization error over 10 different runs in time for filtering and smoothing algorithm}\label{fig:sim1}
\end{figure}
\end{center}

\begin{center}
\begin{figure}[H]
\includegraphics[width=8cm]{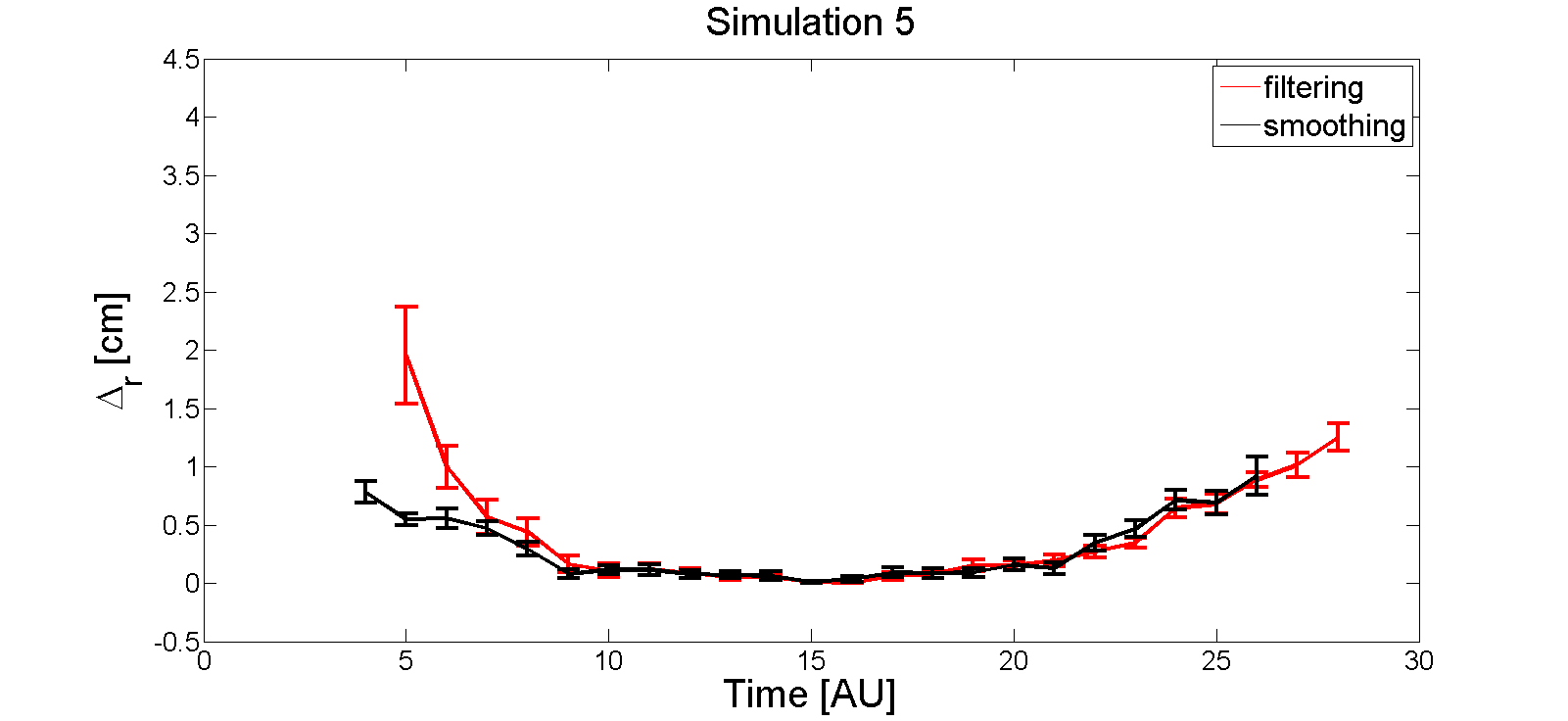}
\includegraphics[width=8cm]{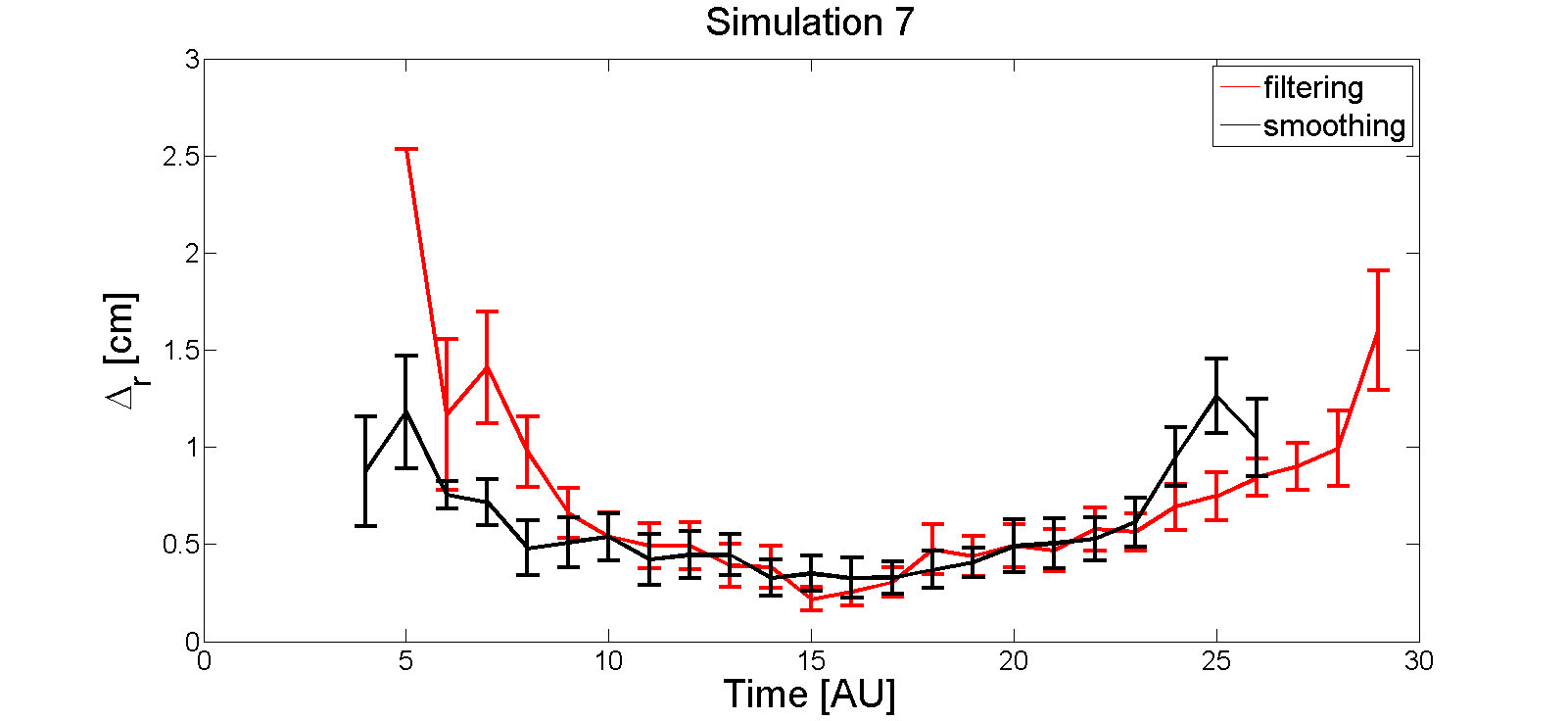}\\ 
\includegraphics[width=8cm]{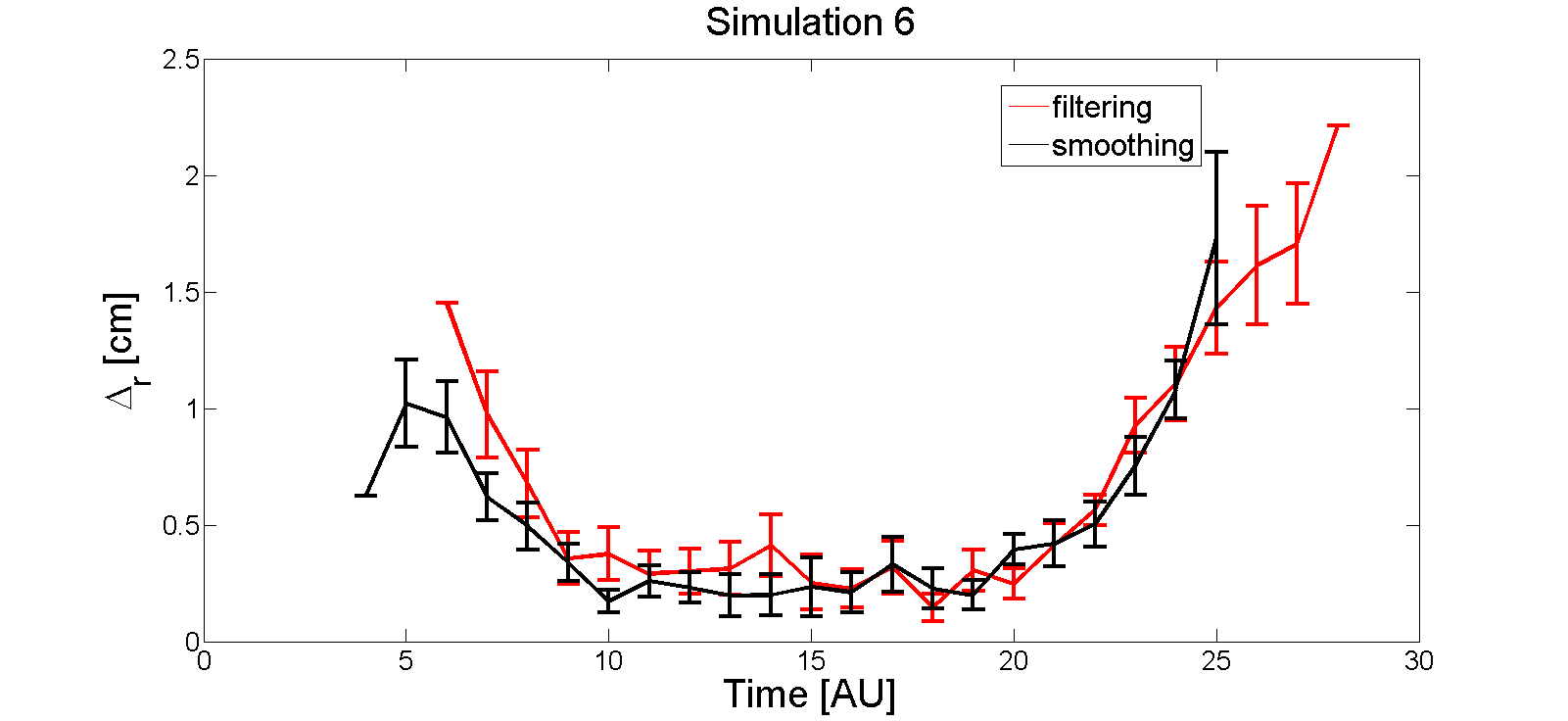}
\includegraphics[width=8cm]{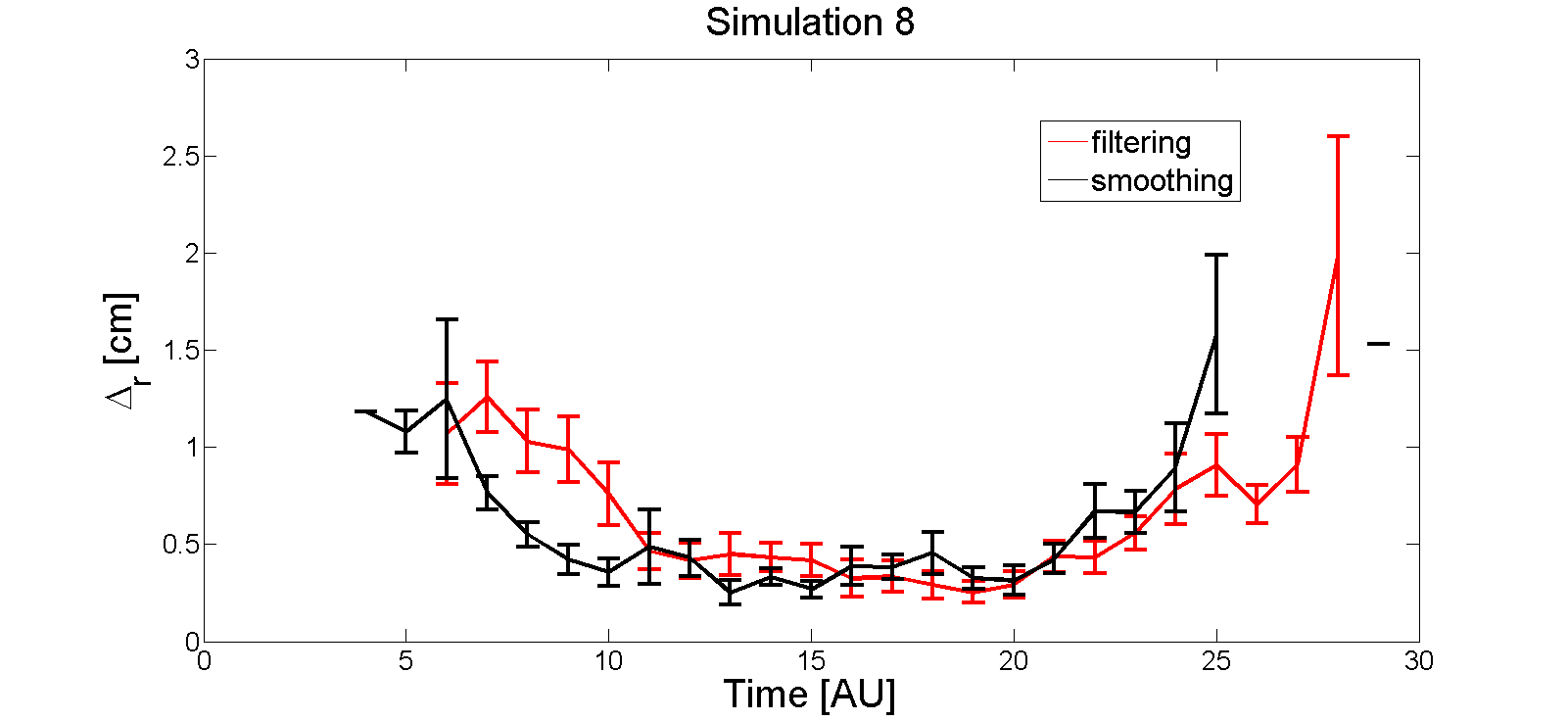}
\caption{Simulation 5--8: Mean localization error over 10 different runs in time for filtering and smoothing algorithm}\label{fig:sim2}
\end{figure}
\end{center}

For each simulation we run our double two--filter smoothing.
We compare the performances of the smoothing algorithm with those of the filtering algorithm, by looking at the localization error, i.e. the distance between the true and the estimated source locations. This is a non-trivial task when the number of estimated dipoles differs from the true one. Following \cite{soetal13, soluar14}, at every time step, we quantify the localization error using a modified version of OSPA \cite{scvovo08} with no penalty for cardinality errors, which are evaluated separately.
Let $(N_t, r_t, {\bf q})$ and $(\hat{N_t}, \hat{r}_t, \hat{{\bf q}} )$ be the true and the estimated dipole configuration respectively; we calculate $\Delta_r$ as 
\begin{equation}
\Delta_r : = \begin{cases}
\min_{\pi \in \Pi_{\hat{N_t}, N_t}}{1 \over \hat{N_t}}\sum_{i=1}^{\hat{N}_t} |\hat{r}_t^{(i)}-{r}_t^{\pi(i)}|&\text{if $\hat{N}_t \leq N_t$}\\

\min_{\pi \in \Pi_{N_t, \hat{N_t}}}{1 \over \hat{N_t}}\sum_{i=1}^{{N}_t} |\hat{r}_t^{\pi(i)}-{r}_t^{(i)}|&\text{if $\hat{N}_t > N_t$}
\end{cases}
\end{equation}
where $\Pi_{k,l}$ is the set of all permutations of $k$ elements drawn from $l$ elements.
In Figure \ref{fig:sim1} and \ref{fig:sim2} we compare the localization error in filtering (red line) and smoothing (black line), averaged over the twenty simulations of each group. The error bars are calculated as the ratio between the standard deviations of the reconstructions and the number of runs in which we have a reconstruction.

A few observations are in order. 
First, for all Simulations the localization error of the smoothing algorithm is systematically lower than that of the filtering algorithm in the first half of the time window. This was indeed expected, as the filtering distribution only embodies information contained in the previous time points, while the smoothing distribution uses (in principle) the whole time series.
A second observation is that in the second half the relative performances of the smoothing algorithm tend to get worse: either the localization error remains comparable to that obtained by the filtering algorithm, or it even becomes larger, particularly for bell--shaped dipoles, where the signal--to--noise ratio decreases at the end of the time--series. This behaviour is most likely explained as a failure in the approximation of the backward information filter, that compromises the approximation of the smoothing distribution. Indeed, the approximation of the backward information filter is certainly worse than that of the forward filter: this fact is confirmed by the asymmetry of the localization error bars of the smoothing algorithm.
A third comment is that there is a manifest difference between the localization errors plotted in Figure \ref{fig:sim1}, where both error bars start from the first time point, and those in Figure \ref{fig:sim2}, where the filtering line starts around $t=5$, while the smoothing line starts from the first time point. This is due to the fact that for bell--shaped dipole moments (those in Fig. \ref{fig:sim2}) the signal strength is relatively low at the beginning, and the filtering algorithm is not capable of estimating the dipole source until around $t=5$; on the contrary, the smoothing algorithm exploits information from the subsequent time points and is therefore able to localize the source.

\section{Experimental Data}
\label{Sec:REA}
We applied the smoothing algorithm to an experimental data set taken from the BESA example database (BESA Gmbh, Munich). Data consist of 32--electrodes EEG recordings from an epileptic subject. 
164 spikes were recorded and averaged, using the peak as trigger (i.e. the point defining $t=0$).
The EEG signals were sampled at 320 Hz and filtered with a Butterworth forward high-pass filter with cut-off frequency of 5 Hz.
The head model is a three--layer model including the brain, the skull and the scalp; while there is no cerebro--spinal fluid (CSF) in the model, the effect of the CSF is partly accounted for by assuming an anisotropic skull conductivity. The tangential conductivity within the skull is modeled to be 3 times larger than the radial conductivity across the skull. The bone conductivities are adjusted to the age of the subject, in this case between 8 and 10 years.
\\

\begin{center}
\begin{figure}[H]
\begin{tabular}{cc}
\includegraphics[width=8cm]{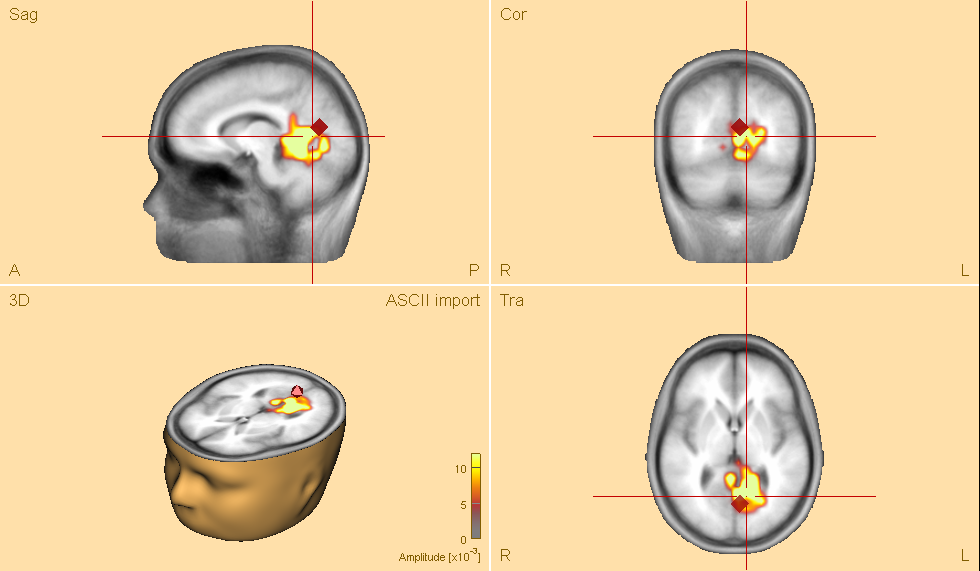} & \\
\includegraphics[width=8cm]{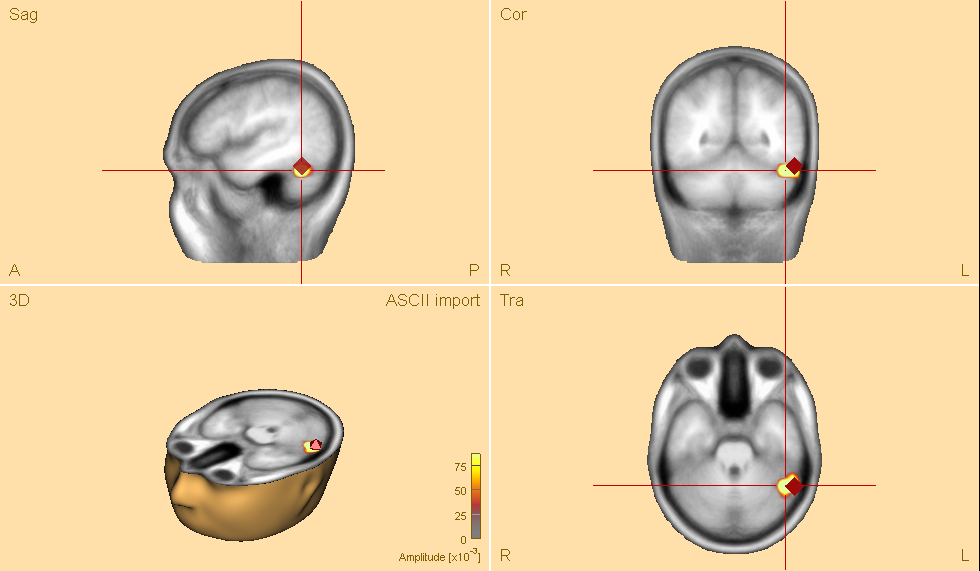} &
\includegraphics[width=8cm]{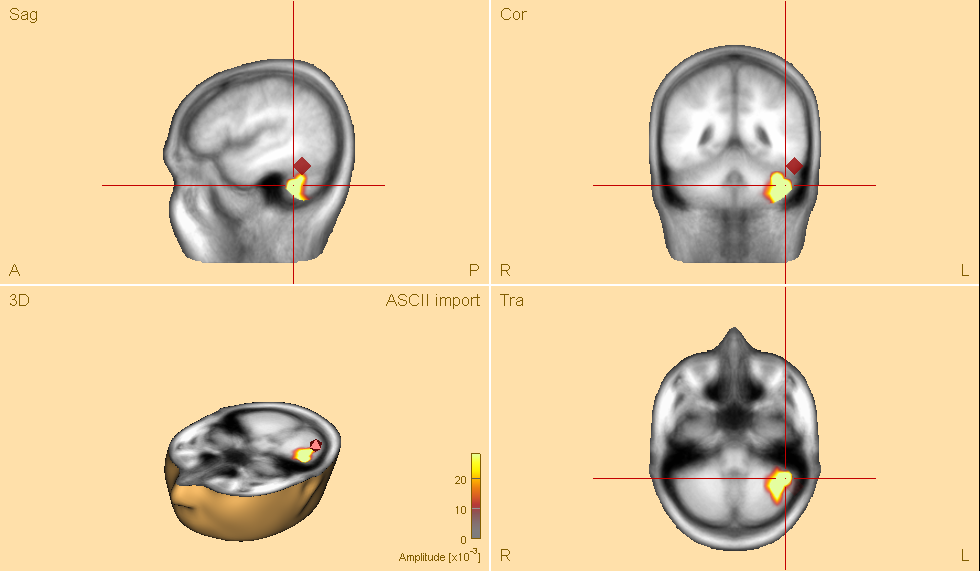}\\ %\vspace{.4cm}
\includegraphics[width=8cm]{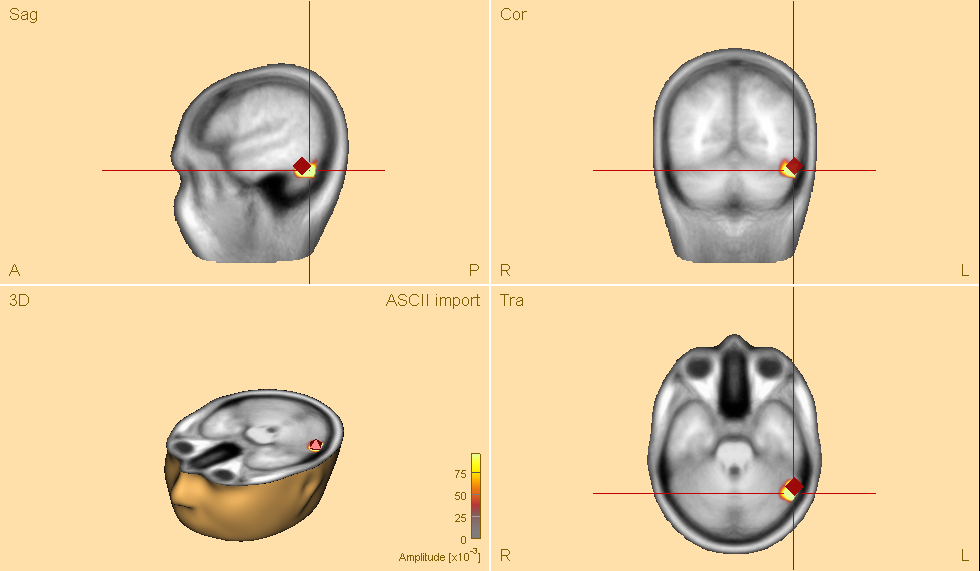} &
\includegraphics[width=8cm]{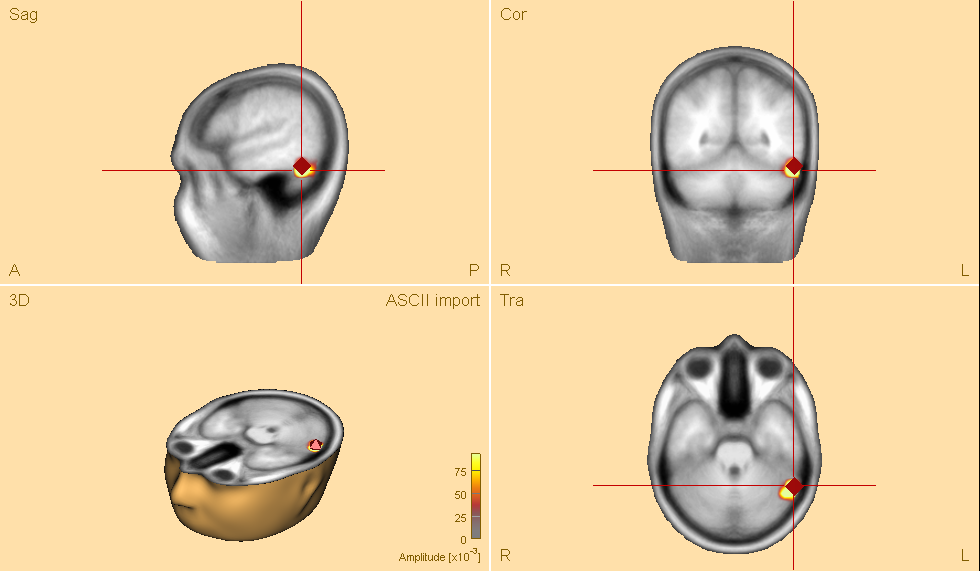}\\ %\vspace{.4cm}
\end{tabular}
\caption{Probability maps of source locations, obtained with the smoothing (left column) and filtering (right column) algorithms at different time points: $t=-40$, $-15$ and 0 ms in the first, second and third row, respectively. }
\label{fig:real_data}
\end{figure}
\end{center}

\noindent
In Figure \ref{fig:real_data} we compare the probability maps obtained by the smoothing and by the filtering algorithm for three selected time points, superimposed on the subject's brain and shown as color maps. For validation, we also show a red diamond corresponding to the location of the dipole estimated by a user--supervised dipole fitting algorithm, applied by an expert user. \\
\noindent
The first time point of interest is $t=-40$ ms, corresponding to the onset of the spike; here the filtering algorithm does not find any source, while the smoothing algorithm obtains a fairly widespread distribution whose support includes the location of the dipole estimated by the expert user. \\
\noindent
The second time point, $t=-15$ ms, corresponds to the propagation of activity from the first location to the actual peak location; here, the probability map obtained by the smoothing algorithm is rather peaked around the location of the dipole estimated by the expert user, while the probability distribution provided by the filtering algorithm is more widespread and peaks a couple of centimeters off the manually estimated dipole.\\
\noindent
Finally, at $t=0$, corresponding to the peak of the spike, the smoothing and the filtering algorithm provide almost identical maps, nicely coherent with the location of the estimated dipole.

\begin{center}
\begin{figure}[H]
\includegraphics[width=14cm]{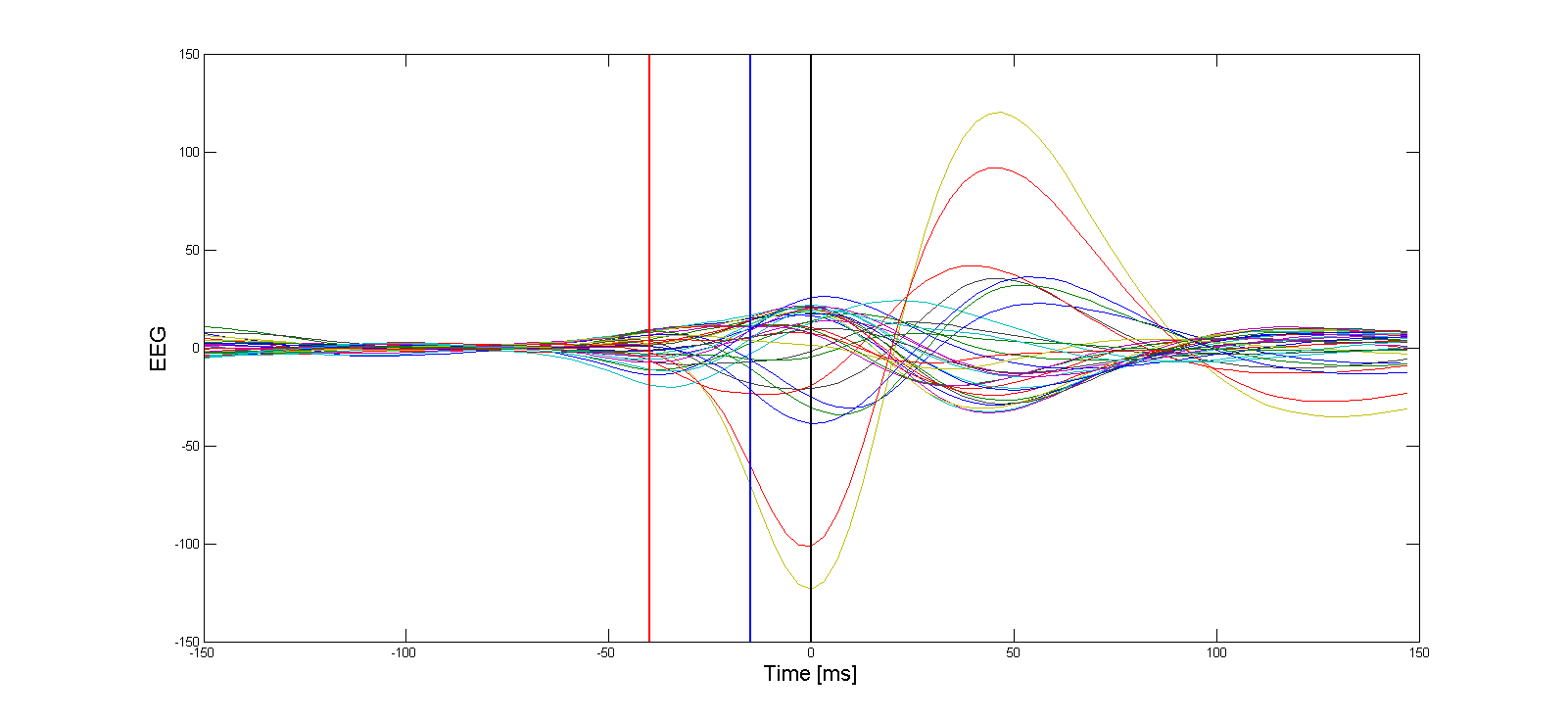}\\ 
\caption{Experimental EEG data taken from the BESA example dataset. The three vertical lines correspond to the time points analyzed in Fig. \ref{fig:real_data}.}
\end{figure}
\end{center}

\section{Discussion}
\label{Sec:DIS}
In this paper we have described a smoothing algorithm, based on the well--known two--filter smoother, for solving the Bayesian inverse M/EEG problem. Our approach obtains two separate approximations of the smoothing distribution, based on the samples of the forward filter and those of the backward information filter, and then selects either of the two, based on the marginal likelihood of the underlying filters. \\
First, we have validated our method by means of eight different synthetic experimental setups, that included dipoles with fixed and moving locations, dipoles with fixed and varying dipole moments, and either one or two simultaneous sources. For each setup, we performed 30 different simulations, by randomly drawing source locations, thus implicitly modifying the SNR of the data. We have confirmed that our approach improves substantially the localization of the sources at their onset, when compared to the filtering, in all the cases under consideration. In particular, for sources with time--varying intensity (with bell--shaped time courses) the smoothing is capable of localizing the source some time points before the filtering; for sources with constant strength, the two algorithms find the source at the same time point but the smoothing has lower average localization error. On the other hand, we have noticed that the approximation of the backward information filter is not as good as that of the forward filter; further work might be devoted to devising better auxiliary distributions to improve the approximation of the backward information filter.\\
Then, we have tested our smoothing algorithm on an experimental data set taken from the BESA example data set,
and we have shown selected time points, that have been chosen as they best represent how the filtering and the smoothing distributions differ at some time points and tend to coincide at others. We have confirmed the superiority of the smoothing algorithm in estimating the onset of the neural sources, by visually comparing the probability maps provided by the smoothing and the filtering algorithm with the source location obtained by an expert user with standard dipole fitting techniques. \\
\noindent
The proposed smoothing algorithm adds to the available tools for source recontruction from M/EEG data, with potentially interesting applications for epilepsy studies where estimating the correct source of epileptic activity, particularly for non-trivial cases such as propagating spikes, is paramount.\\

\noindent
The work described in this article nicely compares to current literature on the M/EEG inverse problem, where the search for temporally smooth solutions is a quite hot research topic. In \cite{ouhago09,gretal12}, in the framework of regularization for distributed source models, the authors propose to use mixed norms (namely, an $L^1$ norm in the spatial domain and an $L^2$ norm in the temporal domain) to incorporate prior knowledge about the continuity of the source time courses. A similar method is presented in \cite{tietal13},
where a functional containing two penalty terms is optimized by means of multivariate penalized regression.
Methods that are conceptually even closer to our smoothing are proposed in \cite{loetal11,laetal12}: here, Bayesian filters for distributed source models are proposed, together with fixed--interval smoothers. Due to the linear/Gaussian model, the authors need only to compute the mean and covariance of the filtering/smoothing distributions, the main difficulty being the size of the state--space, which is large.
To the best of our knowledge, our work is unique in trying to approximate the smoothing distribution for a dynamic set of current dipoles, rather than for a distributed current.\\

\noindent
The algorithm presented in this work is strongly based on the recent literature on sequential Monte Carlo methods; the same literature can provide ideas for further developments. In addition to the already mentioned improvement of the auxiliary distributions, possible future work might include: avoiding the subsampling of the filtering distribution, by means of known strategies for approximating the $\alpha^2$ calculations with a cost of $\alpha \log(\alpha)$ \cite{kletal06}; alternative strategies for sampling the smoothing distribution with linear cost \cite{fewyta10}; finally, exploiting the conditional linearity with respect to the dipole moments, using Rao--Blackwellized smoothing \cite{lietal15}.

\section*{Acknowledgements}

The authors kindly acknowledge the Gruppo Nazionale per il Calcolo Scientifico for financial support. Andre Waelkens, Todor Jordanov and all the BESA staff are kindly acknowledged for their technical and scientific support in the analysis of the experimental data. Finally, we would like to thank the reviewers of the original manuscript for their comments and suggestions, that have improved the quality of this work.

\section{References}
\bibliography{biblio}

\end{document}